\documentclass[12pt]{article}
\usepackage{amsmath}
\usepackage{amsfonts}
\usepackage{amssymb}
\usepackage{a4wide}

\def\W{\mathcal W}
\def\Ev{E^{evil}}

\def\eps{\epsilon}

\def\qed{\hfill \square \ }

\def\C{\mathbb C}
\def\Q{\mathbf Q}
\def\Z{\mathbf Z}

\newtheorem{theorem}{Theorem}[section]
\newtheorem{lemma}[theorem]{Lemma}

\newtheorem{conj}[theorem]{Conjecture}

\title{Irrationality of certain  $p$-adic periods for small $p$}
\author{Frank Calegari\footnote{Supported in part by the American
Institute of Mathematics}}
\begin{document}
\maketitle

\section{Introduction}

Ap\'{e}ry's  proof~\cite{vP} of the irrationality of
$\zeta(3)$ is now over $25$ years old, and it is perhaps surprising that his
methods have not yielded any significantly new results
(although further progress has been made on the irrationality
of zeta values~\cite{R},~\cite{Z}).
Shortly after the initial proof, Beukers produced two elegant
reinterpretations of Ap\'{e}ry's arguments; the first using iterated
integrals and Legendre polynomials~\cite{B1}, and the second using
modular forms~\cite{B2}. It is this second argument that we shall
apply to study the irrationality of certain $p$-adic periods,
in particular, the $p$-adic analogues of $\zeta(3)$ and Catalan's
constant for small $p$.
To relate certain classical
periods to modular forms, Beukers considers various
integrals of holomorphic modular
forms that themselves satisfy certain functional equations (analogous to the
functional equation for the non-holomorphic Eisenstein series of weight two).
That these integrals satisfy functional equations is a consequence
of the theory of Eichler integrals. The periods arise as coefficients
of the associated period polynomials.
In our setting these auxiliary functional equations are replaced by the
notion of 
overconvergent $p$-adic modular forms~\cite{Katz},~\cite{C1},~\cite{C2}. In this guise our
$p$-adic periods will occur as coefficients of overconvergent Eisenstein
series of negative integral weight. These $p$-adic periods are equal
to special values of Kubota--Leopoldt $p$-adic $L$-functions.
Thus $\zeta(3)$ is replaced by $\zeta_p(3) = L_p(3,\mathrm{id})$ and
Catalan's constant 
$$G = L(2,\chi) =  \sum_{n=0}^{\infty}\frac{(-1)^n}{(2n+1)^2}$$
is replaced by $L_p(2,\chi)$, where $\chi$ is
the character of conductor $4$. We shall prove that $\zeta_p(3)$
is irrational for $p=2$ and $3$, and
that $L_p(2,\chi)$ is irrational for $p=2$.

\section{Elementary Remarks and Definitions}

\subsection{Irrationality in $\Q_p$} 

Let $p$ be prime.
Let $\| \cdot \|_p$ denote the $p$-adic absolute value
normalized by $\| p \|_p = 1/p$, and let $| \cdot |$  denote
the Archimedean absolute value.
Given $r = a/b \in  \Q$, how well can $r$ be approximated $p$-adically
by other (distinct) rational integers?

\begin{lemma} Suppose that
$$\left\| \frac{a}{b} - \frac{c}{d} \right\|_p \le \frac{1}{p^n}.$$
Then $\max\{|c|,|d|\} \ge p^n/(|a| + |b|)$.
\label{lemma:approx}
\end{lemma}

\begin{Proof} The inequality above implies that $ad - bc \equiv 0 \mod p^n$.
In particular, it must be the case that $|a d - b c| \ge p^n$, and the
lemma readily follows. $\qed$
\end{Proof}

\medskip

An element $\eta \in \Q_p$ is irrational if it does not lie in $\Q$.
From Lemma~\ref{lemma:approx} we may derive a simple 
criterion for irrationality.

\begin{lemma}[Criterion for $p$-adic irrationality]
Let $p_n/q_n$ be rational numbers with $q_n$ unbounded
and suppose
there exists a $\delta > 0$ such that
$$0 < \left\| \eta - \frac{p_n}{q_n} \right\|_p \le \frac{1}{(\max\{|p_n|,
|q_n|\})^{1 + \delta}}$$
for sufficiently large $n$. Then $\eta$ is irrational.
\end{lemma}

\subsection{Overconvergent Eisenstein Series and $p$-adic Zeta Functions}

Let $B_n$ denote the $n$th Bernoulli number, defined by the identity
$$\frac{x}{2} + \frac{x}{e^x - 1} = \sum_{n=0}^{\infty} \frac{B_n x^n}{n!}.$$
It is a classical result of Euler and Riemann that for non-negative integers $k$,
$$\zeta(1-2k) = - \frac{B_{2k}}{2k}$$
where $\zeta$ is the Riemann zeta function.
Let $\zeta^*_p(s) = (1 - p^{-s}) \zeta(s)$. It is a consequence
of the Kummer congruences that the values $\zeta^*_p(s)$ at
negative odd integers $p$-adically interpolate. In our context
these numbers arise as constant terms of Eisenstein series.
Let $q = e^{2 \pi i \tau}$ and 
suppose $2k \ge 2$ is an even integer. 
Let  $\sigma^*_{2k-1}(n) = \sum_{d|n}^{d \nmid p} d^{2k-1}$. Then
the Eisenstein series: 
$$E^*_{2k}(\tau) = 
\frac{\zeta^*_p(1-2k)}{2} + \sum_{n=1}^{\infty} q^n \sigma^*_{2k-1}(n)$$
is a holomorphic modular form of weight $2k$ for the group $\Gamma_0(p)$.
The form $E^{*}_{2k}$ is related to the more familiar Eisenstein series
of level $\Gamma_0(1)$:
$$E_{2k}(\tau) = \frac{\zeta_p(1-2k)}{2} + \sum_{n=1}^{\infty} q^n \sigma_{2k-1}(n)$$
by the 
relation $E^*_{2k}(\tau) = E_{2k}(\tau)  - p^{2k-1} E_{2k}(p\tau)$.
If $2k$ and $2k'$ are integers congruent modulo $p-1$ and close
$p$-adically, then the $q$-expansions of $E^*_{2k}$ and
$E^*_{2k'}$ are highly congruent. Thus by considering certain
limits of $q$-expansions one may define $p$-adic Eisenstein
series that are $p$-adic modular forms in the sense of 
Serre~\cite{Serre}.
However, one can also make a more precise analytic statement. The 
usual context in which
to view $p$-adic families of $p$-adic eigenforms with finite
slope (non-zero $T_p$-eigenvalue) is Coleman's theory of
overconvergent modular forms~\cite{C1},~\cite{C2} (see also~\cite{Katz}).
This theory
is quite extensive 
and we must be content with the following
(very brief) explanation.   Holomorphic modular forms are
given by sections of certain line bundles over modular curves,
and their $q$-expansions are obtained by evaluating 
at the cusp at infinity. A $p$-adic modular form of
integral weight is a section
of the same sheaf, except considered only over a certain region
of the modular curve defined by points whose corresponding
elliptic curve is ordinary
(more precisely the connected component of this space containing
the cusp at infinity). This region is not open in the sense
of the Zariski topology, but rather in the sense of
rigid analytic geometry (the analytic geometry of the $p$-adic world).
 Finally, an overconvergent modular form of integral weight
is a $p$-adic modular form which extends as a section to some rigid
analytic region of the modular curve strictly containing the
(connected component containing infinity of) the ordinary locus.
It is precisely this overconvergence which will be the $p$-adic
analogue of Beukers' construction whereby a certain analytic
expression $A(z) + \theta B(z)$ converges further for some special
value of $\theta$ than for any other.

\medskip

In general, the weight of a modular form can vary over a rigid
analytic space $\W$ (weight space) of characters
(see~\cite{eigencurve} p. 27).
Let $q = p$ if $p$ is odd and $q = 4$ if $p = 2$. 
 A $\C_p$-point of $\W$
corresponds to a map
$$\Lambda =
\Z_p[[(\Z/q \Z)^{\times} \times (1 + q \Z_p)]]  \simeq
\lim_{\leftarrow} \Z_p[[(\Z/p^n \Z)^{\times}]] \rightarrow \C_p.$$
The classical even integral weights correspond to the characters that
send $a$ with $a \in \lim_{\leftarrow} (\Z/p^n \Z)^{\times}$ to $a^{2k}$.
Let $\W^{+}$ be the part of weight space consisting of characters
$\kappa$ such that $\kappa(-1) = 1$. 
For $n \ge 1$ let
$$\sigma^*_{\kappa}(n):= \sum_{d|n}^{(d,p)=1} \kappa(d) d^{-1}.$$

\begin{theorem} There exists a  function $\zeta_p(\kappa)$ on $\W^{+}$
which is rigid analytic outside $\kappa = 1$ and has a simple pole
at $\kappa = 1$, such that if
$$E_{\kappa} = \frac{\zeta_p(\kappa)}{2} + \sum_{n=1}^{\infty} \sigma^*_{\kappa}(n) q^n$$
Then $E_{\kappa}$ varies rigid analytically over weight space $($away from $\kappa = 1)$
and for each point $\kappa \in \W^{+}$ specializes to an overconvergent eigenform.
Moreover, if $2k \ge 2$ is an even positive integer and
$\kappa$ is the character sending $a$ to $a^{2k}$ 
 then
$E_{\kappa} = E^*_{2k}$.
\end{theorem}

This theorem in this generality is essentially proved in~\cite{C1}.
By continuity, we may recover the value of $\zeta_p(\kappa)$ 
at characters $a \mapsto a^{-2n}$ for positive integers $n$ as
follows:

\begin{lemma} Let $\zeta_p(1+2n)/2$ be the constant term of
$E_{\kappa}$ at the character $a \mapsto a^{-2n}$. Then
$$\zeta_p(1+2n) = \lim_{k \rightarrow n} \zeta(1 - 2k).$$
Where the limit runs over strictly increasing integers $k$
approaching $n$ in $\Z_p$ with the added restriction that
$2k \equiv 2n \mod p-1$.
\end{lemma}

The added restriction comes from the fact that the
characters $a \mapsto a^{2k}$ and $a \mapsto a^{2k'}$ are
close if $2k \equiv 2k' \mod (p-1) p^n$  for large $n$,
as follows from Euler's version of Fermat's little theorem. 
The imposition that $k$ approaches infinity means we
can neglect the Euler factor term, which tends to $1$.

\medskip

\noindent \bf Remark\rm. As with the classical zeta values, not much is known about the
arithmetic nature of
$\zeta_p(1+2n)$. Indeed it could be argued that the situation is worse, as little
is known about the following conjecture, even for $n=1$:

\begin{conj} For all integers $n > 0$ and primes $p$, the values
$\zeta_p(1+2n) \ne 0$. 
\end{conj}

It was the (failed) attempt to prove this conjecture for $n=1$ that lead
to this paper, the idea  being that if one can prove that $\zeta_p(3)$ is
irrational then one has also shown it is non-zero!

\medskip

Let $2k$ be a positive even integer $\ge 4$. If $E_{2k}$ the 
the classical Eisenstein series of weight $2k$ and level $\Gamma_0(p)$
then associated to $E_{2k}$ there is another Eisenstein series,
the \emph{evil twin} $\Ev_{2k}$. The Eisenstein series $\Ev_{2k}$ is
classical, cuspidal at the cusp at $\infty$ and has slope
$2k-1$. Explicitly, the evil twin is given by the following
formula:
$$\Ev_{2k}(\tau) = E_{2k}(\tau) - E_{2k}(p \tau).$$
Let us consider the specialization $E_{-2k}$ of our
Eisenstein family. Let $\theta$ be the operator on
$q$ expansions that acts as $q \cdot d/dq = (2 \pi i)^{-1} d/d\tau$. Then
one can easily compute that
$$\theta^{2k+1} E_{-2k} = \Ev_{2k+2}.$$
Thus we can ``almost'' reconstruct $E_{-2k}$ from
$\Ev_{2k+2}$ by considering
$$E'_{-2k}(\tau) =  (2 \pi i)^{2k+1} \idotsint E_{2k+2} (d \tau)^{2k+1}
\in \Q[[q]]$$
The function $E'_{-2k}$ is holomorphic in a neighbourhood of the
cusp $i \infty$, but is not modular (although similar modified forms
satisfy  some form of functional equation --- see~\cite{B2}, p. 273).
Actually, for our purposes, we could have introduced the
function $E'_{-2k}$ directly without reference to $\Ev_{2k+2}$.
However, by writing down the connection we stress the  
ties between our method and that of Beukers~\cite{B2}.
Let $\eta \in \C_p$ and consider the expression
$$H = E^*_{2k} (E'_{-2k} + \eta).$$
If $\eta = \zeta_p(1+2k)/2$ then $H$ is equal to
$E^{*}_{2k} E_{-2k}$ and is thus an overconvergent modular
function of weight $0$. For all other $\eta$, however, $H$ is
not overconvergent since otherwise $E^*_{2k}$ would be overconvergent
of weight $0$, which is impossible (this follows from~\cite{Katz}, 4.4). Thus we obtain a strictly
analytic (over $\C_p$) characterization of $\zeta_p(1+2k)$.

Suppose that $X_0(p)$ has genus zero. Then the ordinary
locus of $X_0(p)$ containing the cusp at infinity is a rigid
analytic disc.

Suppose that $z$ is a classical
meromorphic modular form that is a local parameter over
the cusp $i \infty$ and has no poles on the component
of the ordinary locus of $X_0(p)$ containing $\infty$ (so $z$
it is overconvergent).
Viewing $z$ first as a complex analytic function, by the
inverse function theorem we may expand
$$H = \sum_{n=0}^{\infty} (a_n - b_n \eta) z^n.$$
Now considering this sum $p$-adically, we know that
$H$ is overconvergent if and only if $\eta = \zeta_p(1+2k)/2$.
Thus we expect the radius of convergence to jump at this point,
and correspondingly the sequence $a_n/b_n$ to converge
$p$-adically to $\zeta_p(1+2k)/2$. If we can estimate
both the $p$-adic and Archimedean radii of convergence for
various $\eta$ we may be able to apply our criterion
of irreducibility. In the next section we carry this
out in detail for $p=2$.

\section{The irrationality of $\zeta_p(3)$  for
$p=2$ and $p=3$}

\subsection{$p = 2$}

Let $p=2$. Then $X_0(2)$ has genus zero, and is uniformized by
the function
$$f = \frac{\Delta(2 \tau)}{\Delta(\tau)} = q \prod_{n=1}^{\infty} (1 + q^n)^{24}.$$
Moreover we note the following facts about $f$ and $X_0(2)$.
The curve $X_0(2)$ has two cusps, $i \infty$ and $0$.
The value of $f$ at these cusps is equal to $0$ and $\infty$, respectively.
In particular, $f$ is a local uniformizer at the cusp at infinity.
The curve $X_0(2)$ has one elliptic point, at $(1 + i)/2$. The value of
$f$ at this point is equal to $-2^{-6}$. This can be proved by noting
that $f'/f = E^*_2$ and that
$$\frac{E^6_2}{\Delta} = \frac{(1 + 2^6 f)^3}{f}.$$
The ordinary locus of (the $2$-adic rigid analytic curve)
$X_0(2)$ has two components, given by 
$$\|f\|_2 \le 1, \qquad
\|f\|_2 \ge 2^{12}.$$
The Fricke involution permutes these two
spaces, as it sends $2^{12} f$ to $1/f$.

\medskip

Consider the series
$$H = E^*_{2k}(E'_{-2k} + \eta) =: \sum_{n=0}^{\infty} (a_n - b_n \eta) f^n.$$
Since  $f = q + \ldots$ it follows that  $b_n \in  \Z$.
Since the $q^n$ coefficient of $E'_{-2k}$ lies in $\Z/n^{2k-1}$ it
also follows that
$[1,2,\ldots,n]^{2k+1} a_n \in \Z$,
where $[1,2,\ldots,n]$ is the greatest common divisor of $1$ up to $n$.

\medskip

Let us  establish the $2$-adic  convergence of this series
for various values of $\eta$.

\begin{lemma} If $p=2$ and
$\eta = \zeta_p(1+2k)/2$ then the radius
of convergence of $H$ is at least $2^{12}$. If $\eta \ne \zeta_p(1+2k)/2$
then the radius of convergence is at most $1$.
\label{lemma:2adic}
\end{lemma}

\begin{Proof} If $\eta = \zeta_p(1+2k)/2$, then $H = E^{*}_{2k} E_{-{2k}}$. Since
$E_{-2k}$ is overconvergent, it extends as a rigid analytic function
somewhere into the supersingular annuli. However, it is also a finite slope eigenform
of level $\Gamma_0(2) = \Gamma_0(p)$,
and such sections  extend far into the supersingular annuli. In particular, by
Theorem~5.2 of Buzzard~\cite{wild}, it extends entirely over the supersingular annuli. Thus
the radius of convergence is at least $2^{12}$. If $\eta \ne \zeta_p(1+2k)/2$,
then $H$ is not overconvergent. Thus it cannot extend into the supersingular
annuli and the radius of convergence is at most $1$. $\qed$
\end{Proof}

From this lemma we may approximate $\zeta_p(1+2k)/2$,
since when $\eta = \zeta_p(1+2k)/2$ the radius of convergence
guaranteed by the previous lemma implies that
$$\|a_n - b_n \eta\|_2 \ll 2^{-(12 - \eps) n}.$$
for any $\eps > 0$ and sufficiently large $n$.
Since $b_n \in \Z$, and since $E^{*}_{2k}$ is not overconvergent of
weight $0$,
it  follows as in the proof of Lemma~\ref{lemma:2adic} that
the $2$-adic valuation of $b_n$ grows slower than any power of $2$.
It follows that (for $\eps > 0$ and $n \gg 0$ as above):
$$\left\|\zeta_2(1+2k) - \frac{2 a_n}{b_n} \right\|_2 \ll 2^{-(12 - \eps) n}.$$

\medskip

Now let us turn to the Archimedean valuations of $a_n$ and $b_n$.
Our arguments here are completely analogous to those of Beukers~\cite{B2}.
By considering $f$ at the cusps and the elliptic points we see that
the radius of convergence of $f$ will be equal to the first branching
value, which occurs at $f((1+i)/2) = -1/2^{6}$. Thus we obtain the estimates:
$$|a_n|,|b_n| \ll  2^{(6+ \eps) n}$$
for all $\eps > 0$ and sufficiently large $n$ (depending on $\eps$). 
The coefficient $a_n$ is not an integer, however. If we write
$a_n = c_n/d_n$ then since $[1,2,\ldots,n]^{2k+1} a_n \in 
\Z$ it follows from the prime number theorem
that (with the usual restrictions on $\eps$ and $n$) that 
$$|d_n| \le [1,2,\ldots,n]^{2k+1} \ll e^{(2k+1+\eps) n}.$$
Consequently, if we write $2 a_n/b_n = p_n/q_n$ where $p_n$ and
$q_n$ are integers, then
$$|p_n| \le |c_n| = |d_n||a_n| \ll 2^{(6 + (2k+1)/\log 2+ \eps) n},
\qquad |q_n| \le |d_n||b_n| \ll   2^{(6 + (2k+1)/\log 2 + \eps) n}.$$
Combining this with our $2$-adic estimates we have proven the following.

\begin{lemma} There exists integers $p_n$, $q_n$ such that
$q_n$ approaches infinity, and such that if
$$\theta = \frac{12 \log 2}{6  \log 2 + 2k + 1}$$
then
$$0 < \left\| \zeta_2(1+2k)  - \frac{p_n}{q_n} \right\|_2 \le \frac{1}{(\max\{|p_n|,
|q_n|\})^{\theta - \eps}}$$
for sufficiently large $n$. 
\end{lemma}

\begin{Proof} This follows from the estimates we have proven so far,
it sufficing to prove that 
 $a_n - \eta b_n \ne 0$ for sufficiently large $n$. Assume
otherwise. Then $H$ is a polynomial in $f$. In particular
$\zeta_p(1+2k) \in \Q$, $H$ has coefficients in $\Q$ and
$E_{-2k}$ is a classical meromorphic eigenform of weight $-2k$.
From the $q$-expansion we may determine that
$E_{-2k}$ has no poles away from the cusps, and a pole of
order at most $1$ at $\tau = 0$. It follows that $H$ is linear
in $f$, which contradicts the fact that $a_2/b_2 \ne a_3/b_3$.
$\qed$
\end{Proof}

As a corollary of this, we prove:

\begin{theorem} If $p=2$ then $\zeta_p(3) \notin \Q$.
\end{theorem}

\begin{Proof}
If $k=1$ then $\theta = 1.1618804316 \ldots > 1$.
Thus we may apply our criterion for irrationality. $\qed$
\end{Proof}

\medskip

If $k=2$ then $\theta =  0.9081638111 < 1$ so we cannot
establish irrationality of $\zeta_2(5)$ (nor indeed $\zeta_2(1+2k)$ for any
other $k \ge 2$).

\medskip

The first few $a_n$ and $b_n$ are given as follows:
$$a_n: 0, 1, 1, -8072/27, 160841/9, -1088512616/1125,
-1088512616/1125 \ldots$$
$$b_n: 1, 24, -552, 19392,  -810024, 37210944, -1815620160 \ldots $$
They are the $2$-adic analogue of Ap\'{e}ry's
sequences
$\{a_n,b_n\}$. 

\medskip

\subsection{$p = 3$}

The same technique can also be applied to other $p$
when $X_0(p)$ has genus zero, where $f$ is chosen
to be $(\Delta(p \tau)/\Delta(\tau))^{1/(p-1)}$.
In this manner we prove the following:

\begin{theorem} If $p=3$ then $\zeta_p(3) \notin \Q$.
\end{theorem}

\begin{Proof}
The construction works as for $p=2$. It suffices to
determine the various radii of convergence. 
If $f = (\Delta(p \tau)/\Delta(\tau))^{1/p-1}$ the
components of the ordinary locus are given
by $\|f\|_3 \le 1$ and $\|f\|_3 \ge 3^6$.
The curve $X_0(3)$ has two cusps $0$ and $i \infty$
at which $f$ has a pole and zero respectively.
There is one elliptic point at $1/2 + \sqrt{-3}/6$,
at which point $f$ takes the value $-3^{-3}$.
Thus one
finds that
$$\theta =\frac{6}{3 + 3/\log 3} = 1.0469892839 \ldots > 1,$$
and thus by the criterion of irrationality we are done.
$\qed$.
\end{Proof}

Although the proof succeeds for $p=3$, it fails for the
other primes where $X_0(p)$ has genus zero ($p=5$, $7$
and $13$).
For example, for $p=5$ we find that
$$\theta = \frac{3}{3/2 + 3/\log 5} = 0.8917942081 < 1.$$

\section{The $2$-adic Catalan's Constant}

There can be no $p$-adic analogue of Ap\'{e}ry's results
for $\zeta(2)$, since $\zeta_p(2) = 0$ for all $p$. However,
we may still study other $p$-adic $L$-values, in particular
the analogue of Catalan's constant:
$$G:= \frac{1}{1^2} - \frac{1}{3^2} + \frac{1}{5^2} - \frac{1}{7^2}
+ \ldots$$
In this context we should study Eisenstein series of \emph{odd}
weight and non-trivial character. Let $p=2$.  Let
$\chi$ be the character of conductor $4$. Then
$L(2,\chi) = G$, whilst for non-negative $k$,
$$L(-2k,\chi) = \frac{E_{2k}}{2}$$
where $E_{2k}$ is the $2k$th Euler number. Moreover, there
exists for each odd positive integer an Eisenstein series
$F_{2k+1} \in S_{2k+1}(\Gamma_1(4),\chi)$ given by the following
$q$-expansion:
$$F_{2k+1} = \frac{L(-2k,\chi)}{2} + \sum_{n=1}^{\infty} q^n
\left( \sum_{d|n} d^{2k} \chi(d) \right)
= \frac{L(-2k,\chi)}{2} + \sum_{n=0}^{\infty} \frac{q^{(2n+1)} (-1)^n (2n+1)^{2k}}
{(1 - q^{2n+1})}.$$
Indeed, $F_{2k+1}$ is nothing but the $2$-adic specialization of
$E_{\kappa}$ to weights of the form $\chi \cdot (a \mapsto a^{2k+1})$,
which are even since $\chi(-1) = -1$.
Since
$$F_{-1} = \frac{L_2(2,\chi)}{2} + \sum_{n=0}^{\infty}  \frac{q^{(2n+1)} (-1)^n} 
{(2n+1)^2 (1 - q^{2n+1})}.$$
We shall consider the function
$$H = F_{1}(F'_{-1} + \eta),$$
where $F'_{-1}$ is holomorphic on the complex upper half plane
and formally given by the $q$-expansion
$$F'_{-1} = \sum_{n=0}^{\infty}  \frac{q^{(2n+1)} (-1)^n}
{(2n+1)^2 (1 - q^{2n+1})}.$$
If $\eta = L_2(2,\chi)/2$ then $F'_{-1} + \eta = F_{-1}$ is overconvergent
(and thus so is $H$), and
otherwise it is not. 
Our arguments now proceed in a very similar manner to those of $\zeta_p(3)$.

\medskip

The curve $X_1(4)$ has genus zero, no elliptic points, and three
cusps corresponding to $i \infty$, $1/2$ and $0$. A uniformizer is
given by
the function
$$z = \left(\frac{\Delta(4 \tau)}{\Delta(\tau)}\right)^{1/3}
=  q \prod_{n=1}^{\infty} (1 + q^n)^8 (1 + q^{2n})^8.$$
The function $z$ vanishes at the cusp at infinity, has a pole
at the cusp at $0$, and equals $-2^{-4}$ at the cusp $1/2$.
The Fricke involution sends $2^8 z$ to $1/z$, and the component
of the ordinary
locus containing infinity is $\|z\|_2 \le 1$.
If we write
$$H = \sum_{n=0}^{\infty} (a_n - b_n \eta) z^n$$
the singularity obtained by inverting $z$ occurs at
$z = -2^{-4}$ and thus
we obtain the 
Archimedean estimates
$$|a_n|, |b_n| \ll 2^{(4 + \eps) n}$$
for all $\eps > 0$ and $n \gg 0$ depending on $\eps$. Furthermore
it is clear that
$$b_n \in \Z, \qquad [1,2,\ldots,n]^2 a_n \in \Z$$\
Thus if $2a_n/b_n = p_n/q_n$ one obtains the estimates
$$|p_n|, |q_n| \ll 2^{(4 + 2/\log 2 + \eps) n}.$$

\medskip

We now study the $2$-adic radii of convergence for various $\eta$.

\begin{lemma} If
$\eta = L_2(2,\chi)/2$ then the radius
of convergence of $H$ is at least $2^{8}$. If $\eta \ne L_2(2,\chi)/2$
then the radius of convergence is at most $1$.
\label{lemma:2adic2}
\end{lemma}

\begin{Proof} If $\eta \ne L_2(2,\chi)/2$ then $H$ is not
overconvergent so the radius of convergence is at most $1$.
Suppose that $\eta = L_2(2,\chi)/2$. Then $F_{-1} = F'_{-1} + \eta$
is an overconvergent finite slope eigenform of level $\Gamma_1(4)$.
Once more we may appeal to the convergence results of Buzzard~\cite{wild},
in particular, Corollary 6.2, to conclude that $F_{-1}$ extends
not only over all of the supersingular region but everywhere
over of $X_1(4)$ except (possibly)  the component of the
ordinary locus containing $0$, which is $\|z\|_2 \ge 2^8$.
$\qed$
\end{Proof}

\medskip

Combining these results, we conclude the following:
\begin{theorem}
There exists integers $p_n$, $q_n$ such that
$q_n$ approaches infinity, and such that if
$$\theta = \frac{8 \log 2}{4  \log 2 + 2} = 1.1618804316 \ldots  > 1$$
then
$$0 < \left\| L_2(2,\chi) - \frac{p_n}{q_n} \right\|_2 \le \frac{1}{(\max\{|p_n|,
|q_n|\})^{\theta - \eps}}$$
for sufficiently large $n$.
In particular, $L_2(2,\chi)$ is irrational.
\end{theorem}

We write down the first few terms $a_n$, $b_n$:
$$a_n: 0, 1, -3, 116/9, -331/9, -99116/225, 3133076/225 \ldots$$
$$b_n: -1, -4,28, -272, 3036, -36624, 464368 \ldots$$
Note that
$$2 \cdot \frac{a_6}{b_6} = \frac{783269}{13060350} =
2^{-1} + 1 + 2^2 + 2^3 + 2^5 + 2^6 + 2^7 + 2^9 + 2^{13} + 2^{18} + \ldots$$
already agrees with $L_2(2,\chi)$ to order $2^{34}$.
As remarked in~\cite{KZ}, the generating function  for $b_n$
automatically
satisfies a second order differential equation.
Indeed one finds that $a_n$ and $b_n$ satisfy the 
Ap\'{e}ry like recurrences:
$$(n+1)^2 u_{n+1} = (4 - 32 n^2) u_n - 256 (n-1)^2 u_{n-1}.$$
If $u_1 = -4$ and $u_{2} = 28$ then $u_n = b_n$, whilst 
if $u_1 = 1$ and $u_2 = -3$ then $u_n = a_n$.

\end{document}